\documentclass[12pt]{amsart}
\usepackage{amsmath}
\usepackage{mathtools}

\usepackage{amsfonts}
\usepackage{tabularx}
\usepackage{amsthm}
\usepackage{autobreak}
\usepackage[english]{babel}
\usepackage{bold-extra}
\usepackage{hyperref}
\usepackage{makecell}
\newtheorem{defi}{Definition}
\newtheorem*{notations}{Notations and conventions}
\newtheorem*{thmA}{Theorem~A}
\newtheorem*{thmB}{Theorem~B}
\newtheorem{prop}[defi]{Proposition}
\newtheorem{cor}[defi]{Corollary}

\newtheorem{lem}[defi]{Lemma}

\newtheorem{thm}[defi]{Theorem}
\newtheorem{question}[defi]{Question}

\newtheorem{rem}[defi]{Remark}
\numberwithin{defi}{section}
\newcommand{\thistheoremname}{}
\newtheorem*{genericthm*}{\thistheoremname}
\newenvironment{namedthm*}[1]
{\renewcommand{\thistheoremname}{#1}%
	\begin{genericthm*}}
	{\end{genericthm*}}
\allowdisplaybreaks

\def\D{{\mathbb{D}}}

\def\R{{\mathbb{R}}}
\def\C{{\mathbb{C}}}

\def\N{{\mathbb{N}}}

\def\T{{\mathbb{T}}}
\def\UU{{\mathcal{U}}}

\def\GG{{\mathcal{G}}}
\def\DD{{\mathcal{D}}}

\def\HH{{\mathcal{H}}}

\usepackage{multicol}
\usepackage[T1]{fontenc}
\usepackage[utf8]{inputenc}
\usepackage[a4paper]{geometry}
\usepackage{amssymb,amsmath}
\usepackage{version}
\usepackage{mathtools} 
\usepackage{enumerate}
\usepackage{dsfont}
\usepackage{mathdots}
\makeatletter

\usepackage{latexsym}
\usepackage{color}
\usepackage[normalem]{ulem}

\definecolor{skyblue}{rgb}{0.0,0.55,0.85}

\definecolor{darkgreen}{rgb}{0,0.38,0}

\definecolor{mybrown}{rgb}{0.59,0.29,0.0}

\begin{document}

\title[Generic divergence of Taylor coefficients under composition]{Generic divergence of Taylor coefficients under composition in weighted Beurling--Sobolev spaces}

\author{St\'{e}phane Charpentier}
\address{Aix-Marseille Univ, CNRS, I2M, Marseille, France}
\email{stephane.charpentier.1@univ-amu.fr}

\author{Nicolas Espoullier}
\address{Aix-Marseille Univ, CNRS, I2M, Marseille, France}
\email{nicolas.espoullier@etu.univ-amu.fr}

\author{Rachid Zarouf}
\address{Aix Marseille Univ, Lab ADEF, Campus Univ St Jerome, 52 Ave Escadrille Normandie, F-13013 Marseille, France}
\address{Aix-Marseille Univ., Univ. de Toulon, CNRS, CPT, Marseille, France}

\email{rachid.zarouf@univ-amu.fr}

\thanks{The first author is supported by the grant ANR-24-CE40-0892-01 of the French National Research Agency
	ANR (project COMOP). The work of Rachid Zarouf is supported by the pilot center Ampiric, funded by the France 2030 Investment Program operated by the Caisse des D\'ep\^ots.}


\begin{abstract}Let $p\geq 1$ and $w$ be a slowly oscillating weight. Let also $(v_k)_{k\geq 0}$ be a sequence of positive numbers such that $v_k\to \infty$. On the one hand, we prove the existence of a comeagre set $\GG$ in the Beurling--Sobolev space $\ell_a^p(w)$ such that, for every $f\in\mathcal G$ and every non-monomial finite Blaschke product $B$, if we write
	\[
	f\circ B =\sum_{k\geq 0}a_kz^k,
	\]
	then
	\[
	\sum_{k\geq 1} |a_kw_k|^{2-\frac{v_k}{\log(k+1)}}=\infty.
	\]
	In particular, this shows that, for a generic $f \in \ell_a^1(w)$, $f\circ B$ fails to belong to any space $\ell_a^p(w)$ with $p<2$. On the other hand, we obtain that, if we assume $\log(v_k)=o(\log(k))$ as $k\to \infty$, then, for every $p>2$ and every non-monomial finite Blaschke product $B$, there exists a comeagre set $\UU_B$ in $\ell_a^p(w)$ such that, for every $f\in \UU_B$, writing $f\circ B=\sum_{k\ge0}b_kz^k$, one has
	\[
	\limsup_{k\to\infty}|b_k|w_kv_k^{-1}=\infty.
	\]
	These strengthen and extend a result due to Hal\'asz.
\end{abstract}

\maketitle

\section{Introduction}

A sequence $w=(w_k)_{k\geq 0}$ of positive real numbers is called a weight if it satisfies
\begin{equation}\label{condition-0-weight}
\liminf_{k\to \infty}w_k^{1/k}\geq 1.
\end{equation}
Letting $\D:=\{z\in \C:\, |z|<1\}$ be the open unit disk of the complex plane, we write $H(\D)$ for the Fr\'echet space of all functions analytic in $\D$.  Given $p\in[1,\infty)$ and  a weight $w$, we denote by $\ell_{a}^{p}(w)$ the weighted Beurling--Sobolev space consisting of all functions $f=\sum_{k\geq 0}a_kz^k$ in $H(\D)$ such that
\[
\Vert f\Vert_{\ell_a^{p}(w)}^{p}:=\sum_{k\geq 0}|a_k|^{p}w_k^p<\infty.
\]
For $p=\infty$ we write $\ell_a^{\infty}(w)$ for the space of all $f\in H(\D)$ for which
\[
\|f\|_{\ell_a^{\infty}(w)}:=\sup_{k\geq 0}|a_k|w_k<\infty.
\]
By \eqref{condition-0-weight},  each space $\ell_a^p(w)$, $p\in [1,\infty]$, embeds continuously into $H(\D)$. Given two weights $w$ and $\beta$, we say that $w$ and $\beta$ are comparable, and we write $w\asymp \beta$, provided there exist two constants $c,C>0$ such that
\[
c \beta _n \leq w_n \leq C \beta_n,
\]
for every $n\in \N$. If $w\asymp \beta$, then $\ell_a^p(w)=\ell_a^p(\beta)$ with equivalent norms. When $w\asymp 1$, we write $\ell_a^p:=\ell_a^p(w)$. Depending on the growth of the weight, Beurling--Sobolev spaces form a class of Banach spaces of analytic functions on $\D$ containing the classical ones. For example, $W:=\ell_a^1$ is the Wiener space and $H^2:=\ell_a^2$ is the Hardy space. If $w\asymp (n+1)^{-1/2}$, then $\ell_a^2(w)$ is the Bergman space $A^2$ and, if $w\asymp (n+1)^{1/2}$, then one recovers the Dirichlet space $\DD:=\ell_a^2(w)$.

\medskip

We recall that, given $\mu \in \partial \D$ and a finite sequence $\Lambda=(\lambda_1,\ldots,\lambda_d)\subset \D$, the finite Blaschke product $B$ associated with $\Lambda$ is defined by
\begin{equation}\label{Blaschke-def}
B(z)=\mu\prod_{j=1}^d\frac{\lambda_j -z}{1-\overline{\lambda_j}z},\quad z\in \D.
\end{equation}
The natural number $d$ is called the degree of $B$. When $d=1$, $B$ is called a single Blaschke factor and it is an automorphism of $\D$. We write $b_{\lambda}$ the single Blaschke factor vanishing at $\lambda$. A monomial is any finite Blaschke product associated to $\Lambda=(0,\ldots,0)$. Like every function in $H(\overline{\D})$, every finite Blaschke product belongs to each $\ell_a^p(w)$ space, if one assumes that $w$ is slowly oscillating.

\medskip

In a paper that seems to have received little attention, Gabor Hal\'asz established the following result.

\begin{thm}[Hal\'asz \cite{Halasz1967}]\label{Halaszthm}Let $(v_k)_k$ be a sequence of positive real numbers satisfying $v_k \to \infty$. For every $\lambda\in\D\setminus\{0\}$, there exists a function $f\in \ell_a^1$ such that, if we write $f\circ b_{\lambda} = \sum_{k\geq 0}a_k z^k$, then
	\[
\sum_{k\geq 2} |a_k|^{2-\frac{v_k}{\log(k)}} =\infty.
	\]
\end{thm}
Note that, by choosing $v_k$ such that $v_k/\log(k) \to 0$, Hal\'asz's result guarantees the existence of $f\in \ell_a^1$ such that
\[
f \circ b_\lambda \notin \bigcup_{q\in[1,2)} \ell_a^q.
\]
In other words, calling a symbol an analytic self-map of $\D$, this theorem tells us that the composition operator
\[
C_{\varphi}:f\mapsto f\circ \varphi
\]
does not map $\ell_a^1$ into $\ell_a^p$, for every $p\in[1,2)$, whenever the symbol $\varphi$ is a Blaschke factor $b_\lambda$ with $\lambda \in \D\setminus \{0\}$. This statement became a direct consequence of a more general result, due to Blyudze and Shimorin in \cite{BlyudzeShimorin1993}, who established a characterisation of those symbols $\varphi \in H(\overline{\D})$ such that $C_\varphi$ maps $\ell_a^1$ into $\ell_a^p$, $p\in [1,2)$. In particular, they obtained, as a corollary, that $C_B$ does not map $\ell_a^1$ into $\ell_a^p$, for every non monomial finite Blaschke product $B$ and every $p\in [1,2)$. In the recent paper \cite{CharpentierZarouf}, the authors extended Blyudze--Shimorin's result to the weighted setting, and established a description of all symbols $\varphi\in H(\overline{\D})$ such that $C_\varphi$ maps $\ell_a^p(w)$ to $\ell_a^{\infty}(w)$, $p>2$, under the assumption that the weight is slowly oscillating; see below for the definition of a slowly oscillating weight. In the same paper, we mention that they obtain a complete characterization of those $\varphi\in H(\overline{\D})$, for which $C_\varphi$ maps $\ell_a^p$ into $\ell_a^q$. It appears that, for every non-monomial finite Blaschke product $B$, $C_B$ does not map $\ell_a^p(w)$ into $\ell_a^{\infty}(w)$, whenever $p>2$ and $w$ is slowly oscillating. Let us indicate that Lefèvre \emph{et al.} proved that $C_{b_\lambda}$ does not map $\ell_a^p(w)$ into itself, whenever $\lambda\neq 0$ and $p\neq 2$, whatever the weight $w$ \cite{LLQR-lpa-2025}.

\medskip

The purpose of this note is threefold: first, to improve Hal\'asz's result by proving that it holds for every non-monomial finite Blaschke product and for certain weighted Beurling--Sobolev spaces; second, to establish a similar statement for finite Blaschke product composition operators from $\ell_a^p(w)$ to $\ell_a^{\infty}(w)$; third, to highlight a subsequent \emph{generic} phenomenon. More precisely, let us divide our contributions into the following two results. We recall that a subset of a topological space is $G_{\delta}$ if it is a countable intersection of open sets. Moreover,  a weight $w$ is slowly oscillating if, for every $0<J<1<L<\infty$, there exists $c> 0$ such that, for every $n\in \N$ and every $Jn\leq k \leq Ln$, one has
\[
c^{-1}w_{n}\leq w_k\leq cw_{n}.
\]

\begin{thmA}Let $w$ be a slowly oscillating weight and let $(v_k)_k$ be a sequence of positive real numbers satisfying $v_k\to \infty$. For every $p\geq1$, there exists a dense $G_\delta$ subset $\mathcal G$ of $\ell_a^p(w)$ such that, for every $f\in\mathcal G$ and every non-monomial finite Blaschke product $B$, we have
	\[
\sum_{k\geq 1} |\widehat{f\circ B}(k)w_k|^{2-\frac{v_k}{\log(k+1)}}=\infty.
	\]
\end{thmA}

Our second main result reads as follows.

\begin{thmB}Let $p>2$, let $w$ be a slowly oscillating weight, and let $(v_k)_{k\geq 1}$ be a sequence of positive real numbers satisfying $v_k\to \infty$. Assume that
	\begin{equation*}
		\frac{\log v_k}{\log(k+1)}\longrightarrow 0\quad \text{as }k\to \infty.
	\end{equation*}
	Then, for every non-monomial finite Blaschke product $B$, there exists a dense $G_\delta$-subset $\HH$ of $\ell_a^p(w)$ such that, for every $f\in\HH$,
	\begin{equation*}
		\limsup_{k\to\infty}\frac{|\widehat{f\circ B}(k)|w_k}{v_k}=\infty.
	\end{equation*}
	In particular, for every $f\in\HH$, one has $f\circ B\notin\ell_a^\infty(w)$.
\end{thmB}

Both results are optimal in some sense; for Theorem~B see Remark~\ref{rem-optimal-subpolynomial} and, for Theorem~A, the optimality of the condition on $(v_k)_k$ is already discussed in \cite{Halasz1967}.

As one can easily see using Baire's theorem, if $C_\varphi$ does not map $\ell_a^p(w)$ into $\ell_a^q(w)$, $p,q\in [1,\infty]$, then the set of all functions $f\in \ell_a^p(w)$ such that $f\circ \varphi \notin \ell_a^q(w)$ is comeagre (Proposition~\ref{prop-general-Cphi-Halasz}). Thus, the most striking feature of Theorem~A is, perhaps, that the comeagre set $\mathcal G$ can be chosen independently of the non-monomial finite Blaschke product $B$.  The question whether this can also be achieved in Theorem~B remains open. Additionally, both Theorems~A and~B provide \emph{sharp} quantitative information on the way $|\widehat{f\circ B}(k)|$ may be large when $k\to \infty$.

While Hal\'asz's proof of Theorem~\ref{Halaszthm} was constructive, ours eventually follows from an application of the Baire Category Theorem. Whatever the constructive or the Baire approach, the key observation is that, given $\varphi \in H(\overline{\D})$,
\[
\forall f\in \ell_a^1(w),\, f\circ \varphi \in \ell_a^p(w) \quad \Longleftrightarrow \quad \alpha_\varphi(w,p):=\sup_{n\geq 0}\sum_{k\geq 0}|\widehat{\varphi^n}(k)|^p\frac{w_k^p}{w_n^p} <\infty
\]
and
\[
\forall f\in \ell_a^p(w),\,f\circ \varphi \in \ell_a^\infty(w) \quad \Longleftrightarrow \quad \beta_\varphi(w,p):=\sup_{k\geq 0}\sum_{n\geq 0}|\widehat{\varphi^n}(k)|^{p'}\frac{w_k^{p'}}{w_n^{p'}} <\infty,
\]
where $p'$ is the conjugate exponent of $p$. Quantitative lower and upper estimates of  $\alpha_\varphi(w,p)$ were obtained by Blyudze and Shimorin in \cite{BlyudzeShimorin1993} for every $p\in[1,2)$ and every $\varphi\in H(\overline{\D})$, using difficult tools from Fourier analysis and asymptotics analysis. Similar estimates were recently established for $\beta_\varphi(w,p)$ in \cite{CharpentierZarouf}, making it appear a subtle difficulty in lower estimating $\beta_\varphi(w,p)$. This difficulty is one of the main obstacles to obtaining a version of Theorem~B uniform in $B$. In order to prove Theorem~A, instead of searching for uniform version of Blyudze--Shimorin estimates of $\alpha_B(w,p)$, which can be rather difficult, we will rather make use of a simple trick, already appearing in \cite{CharpentierEspoullierZarouf2026arxiv}, that reduces our problem to uniform quantitative estimates on the $\ell_a^\infty$ norm of $B^n$  (Theorem~\ref{prop-uniform-main-esti-1} below). For the proof of Theorem~B, since we do not need uniform estimates, we can directly base our arguments on the estimates of $\beta_B(w,p)$ obtained in \cite{CharpentierZarouf}.

\medskip

The paper is organized as follows. In Section~\ref{section2}, we establish the required estimates of the Fourier coefficients of powers of $B$, in particular the uniform estimates. Theorem~A is proved in Section~\ref{sec-Halasz} and the proof of Theorem~B is given in Section~\ref{sec-dual-fixed-B}.

\begin{notations}{\rm Throughout the paper, we use the standard notation $O(\cdot)$ and $o(\cdot)$. For nonnegative quantities $A$ and $B$, we write $A\lesssim B$ if there exists a constant $C>0$ such that $A\leq CB$ and $A\asymp B$ if both $A\lesssim B$ and $A\gtrsim B$ hold. The involved constants may depend on some fixed parameters.
	
We also use the convention $0^{\beta}=+\infty$ if $\beta <0$ and $0^0=1$. }
\end{notations}

\section{Fourier coefficient estimates for powers of finite Blaschke products}\label{section2}

First of all, we will make use of the following elementary result, that can be found in \cite{CharpentierZarouf}.

\begin{lem}[Lemma~2.1 in \cite{CharpentierZarouf}]\label{lem-elementary}Let $p\in [1,\infty)$ and let $\varphi \in H(\overline{\D})$ be  a symbol. Then there exist constants $0<J<1<L<\infty$ such that, for every $\gamma \geq 0$, there exists $s\in(0,1)$ such that, for every $n$ large enough,
	\[
	\sum_{k<Jn}\left|\widehat{\varphi^{n}}(k)\right|^{p}k^\gamma+\sum_{k>Ln}\left|\widehat{\varphi^{n}}(k)\right|^{p}k^\gamma \leq s^n
	\]
	and, for every $k$ large enough,
	\[
	\sum_{n<Jk}\left|\widehat{\varphi^{n}}(k)\right|^{p} n^\gamma + \sum_{n>Lk}\left|\widehat{\varphi^{n}}(k)\right|^{p}n^\gamma \leq s^k.
	\]
\end{lem}

We will also repeatedly use the following immediate application of Cauchy's inequalities, which somehow follows from the proof of the previous proposition, but that we isolate for convenience.

\begin{lem}\label{lem-easy-expo-decay}Let $\varphi$ be holomorphic in a neighborhood of  the closed disc $\overline{D(0,R)}$ for some $R>1$, and let $P$ be a polynomial. Then, for every $k\geq 0$,
\[
|\widehat{P\circ \varphi}(k)|\lesssim R^{-k},
\]
where the constant depends on $R$, $P$ and $\varphi$.
\end{lem}

\begin{proof}Let $\varphi$, $R$ and $P$ be as in the statement. By Cauchy's inequalities, we have, for every $k\geq 0$,
\[
|\widehat{P\circ \varphi}(k)| \lesssim \frac{\max_{|z|\leq R} |P\circ \varphi(z)|}{R^k}.
\]
\end{proof}

We turn to the central ingredients of the proofs of the main results. The following estimate is established in \cite{BFZ-Blaschkeproduct-2024}.

\begin{thm}[Theorem~1 in \cite{BFZ-Blaschkeproduct-2024}]\label{Mainesti}Let $B$ be a non-monomial finite Blaschke product. There exists a constant $N_B \geq 3$ such that, for every $n$ large enough,
	\begin{equation}\label{Main-esti-1}
	\sup_{k\ge0}|\widehat{B^{n}}(k)|\lesssim n^{-1/N_B}.
	\end{equation}
\end{thm}
This estimate was shown to be sharp as $n\to\infty$. For the proof of Theorem~A, we will need a uniform version of it. Let us introduce some notations. For $\mu\in \partial \D$, $d\geq 1$ and $\lambda_1,\ldots ,\lambda_d \in \D$, let $B$ be the finite Blaschke product associated with $\lambda_1,\ldots , \lambda_d$,
\[
B(z)=\mu\prod_{j=1}^{d}\frac{z-\lambda_{j}}{\,1-\overline{\lambda_{j}}\,z\,},\qquad z\in\mathbb{D}.
\]
For the matter we are interested in, the parameter $\mu$ plays no role and, from now on, we may and shall assume that $\mu=1$. For $n,k\in \N$, the $k$th Fourier coefficient of $B^n$ are given by
\[
\widehat{B^{n}}(k)=\frac{1}{2\pi}\int_{0}^{2\pi}B^{n}(e^{i\theta})e^{-ik\theta}\,d\theta.
\]
We denote by $\psi_B:\R\to \R$ the continuous argument of the map $\theta\mapsto B(e^{i\theta})$,
\[
B(e^{{\rm i}\theta})=e^{ i\psi_{B}(\theta)}\quad \text{and}\quad \psi_B(\theta+2\pi) = \psi_B(\theta) +2\pi d,\quad \theta \in \R.
\]
One can easily show (see, for example, \cite{BFZ-Blaschkeproduct-2024}) that
\[
\psi_{B}'(\theta)=\Re\!\sum_{j=1}^{d}\frac{1+\overline{\lambda_{j}}e^{i\theta}}{1-\overline{\lambda_{j}}e^{i\theta}},
\]
so that $\psi_B$ a real-analytic non-decreasing function. Moreover, $\psi_B''$ is not identically zero and has finitely many zeros if and only if $B$ is not monomial. Assume that $B$ is non-monomial and let us then denote by $(\xi_{\ell})_{1\leq \ell \leq s}$ the finite sequence of the consecutive zeros of $\psi_B''$ in $[0,2\pi)$, with respective multiplicities $(N_{\ell}-2)_{1\leq \ell \leq s}$, $N_{\ell}\ge3$. This means that 
\[
0\leq\xi_{1}<\xi_{2}<\ldots<\xi_{s}< 2\pi
\]
and that for any $\ell=1,\ldots,s$ 
\[
\psi''_{B}(\xi_{\ell})=\ldots=\psi_{B}^{(N_{\ell}-1)}(\xi_{\ell})=0,\qquad\psi_{B}^{(N_{\ell})}(\xi_{\ell})\not=0.
\]
We set
\begin{equation}\label{defi_N_B}
	N_{B}=\max_{1\le\ell\le s}N_{\ell}.
\end{equation}

\medskip

Let us now focus on a refinement of Theorem~\ref{Mainesti}~\eqref{Main-esti-1} (and Lemma~\ref{lem-elementary}) that will make explicit a certain uniformity over suitable families of finite Blaschke products, and that will be needed for the proof of Theorem~A. We begin by introducing some additional notation.

\medskip

For $d\ge1$ and $0<\eta<r<1$, set
\begin{equation}\label{defiBdr}
\mathcal B_{d,\eta,r}:=
\left\{
B(z)=\mu\prod_{j=1}^d\frac{z-\lambda_j}{1-\overline{\lambda_j}z}:\,|\mu|=1,\ \eta \leq \max_j|\lambda_j|\leq r 
\right\}
\end{equation}
and
\begin{equation}\label{defiKdr}
K_{d,\eta,r}:=\left\{\Lambda=(\mu,\lambda_1,\ldots,\lambda_d)\in\T\times\D^d:\,\eta \leq \max_j|\lambda_j|\leq r \right\}.
\end{equation}
For $\Lambda=(\mu,\lambda_1,\ldots,\lambda_d)\in K_{d,\eta,r}$, denote by $B_{\Lambda}$ the corresponding Blaschke product in \eqref{defiBdr}.
Note that, for every integer $q\ge1$, the map $(\Lambda,\theta)\longmapsto \psi_{B_\Lambda}^{(q)}(\theta)$ is continuous on the compact set $K_{d,\eta,r}\times[0,2\pi]$. Indeed this can be checked from the explicit formula
\[
\psi'_{B_\Lambda}(\theta)
=\Re\sum_{j=1}^d\frac{1+\overline{\lambda_j} e^{i\theta}}
{1-\overline{\lambda_j} e^{i\theta}},
\]
and its derivatives.

Our refinement of Theorem~\ref{Mainesti}~\eqref{Main-esti-1} reads as follows.

\begin{thm}\label{prop-uniform-main-esti-1}
	For every $d\ge1$ and every $0<\eta<r<1$, there exist $N_{d,\eta,r}\ge3$
	and $C_{d,\eta,r}>0$ such that
	\[
	\sup_{B\in\mathcal B_{d,\eta,r}}\sup_{k\ge0}|\widehat{B^n}(k)|
	\le C_{d,\eta,r} n^{-1/N_{d,\eta,r}},\quad n\ge1.
	\]
	Moreover, there exist $0<J<1<L$, $q\in (0,1)$ and $n_0\in \N$, each of these constants depending only on $d,\eta,r$, such that for every $n\geq n_0$,
	\[
	\sup_{B\in \mathcal{B}_{d,\eta,r}}\sum_{k<Jn}|\widehat{B^n}(k)|^2\leq
	q^n\quad
	\text{and}
	\quad
	\sup_{B\in \mathcal{B}_{d,\eta,r}}\sum_{k>Ln}|\widehat{B^n}(k)|^2\leq q^n.
	\]
\end{thm}

The proof of Theorem~\ref{prop-uniform-main-esti-1} relies on the following compactness lemma, which isolates the key estimate that holds uniformly for finite Blaschke products in $\mathcal B_{d,\eta,r}$.

\begin{lem}\label{lem:uniform-finite-type}
	For every $d\ge1$ and every $0<\eta<r<1$, there exist an integer $N_{d,\eta,r}\ge3$
	and a constant $c_{d,\eta,r}>0$ such that, for every $B\in\mathcal B_{d,\eta,r}$,
	\[
	\max_{2\le j\le N_{d,\eta,r}}|\psi_B^{(j)}(\theta)|\ge c_{d,\eta,r},
	\quad \theta\in[0,2\pi].
	\]
\end{lem}

\begin{proof}
	Assume the contrary. Then, for every integer $N\geq 2$ and every constant $c>0$, there exist $B\in\mathcal B_{d,\eta,r}$ and $\theta\in[0,2\pi]$ such that
	\[
	\max_{2\le j\le N}|\psi_B^{(j)}(\theta)|<c.
	\]
	Applying this statement with $N=m$ and $c=1/m$, we obtain, for every integer $m\ge2$,
	a Blaschke product $B_m=B_{\Lambda^{(m)}}\in\mathcal B_{d,\eta,r}$ and a point
	$\theta_m\in[0,2\pi]$ such that
	\[
	\max_{2\le j\le m}|\psi_{B_m}^{(j)}(\theta_m)|<\frac1m.
	\]
	By compactness of $K_{d,\eta,r}\times [0,2\pi]$, after extraction, we may assume that
	$\Lambda^{(m)}\to\Lambda^{(\infty)}\in K_{d,\eta,r}$ and $\theta_m\to\theta_\infty\in[0,2\pi]$.
	Set $B_\infty:=B_{\Lambda^{(\infty)}}$. Since, for every integer $q\ge1$, the map
	$(\Lambda,\theta)\mapsto \psi_{B_\Lambda}^{(q)}(\theta)$ is continuous on
	$K_{d,\eta,r}\times [0,2\pi]$, passing to the limit gives
	\[
	\psi_{B_\infty}^{(j)}(\theta_\infty)=0,
	\qquad j\ge2.
	\]
	Since 	$\psi_{B_\infty}''$ is real analytic on $[0,2\pi]$, it follows that $\psi_{B_\infty}''\equiv0$, hence $B_\infty$ is monomial. This
	contradicts the fact that $B_\infty\in\mathcal B_{d,\eta,r}$, because, by definition of
	$\mathcal B_{d,\eta,r}$, at least one zero of $B_\infty$ has modulus at least $\eta>0$.
\end{proof}

The key ingredient in the proof of Theorem~\ref{prop-uniform-main-esti-1}  is  the following van der Corput estimate, that we recall below in the form stated by Stein \cite{SteinBook1993}.

\begin{lem}[\cite{SteinBook1993}, Proposition 2, page 332]\label{lem:vdc} Let $j\geq 1$ and let
	$\Phi:(a,b)\to \R$ be a $C^{j}$ function. Assume that, for every $x\in (a,b)$,
	\[
	|\Phi^{(j)}(x)|\geq 1.
	\]
	Then, for every $\nu\in \R$, the estimate
	\[
	\Big|\int_{a}^{b}e^{i\nu\Phi(x)}\,dx\Big|\;\le\;(5\cdot2^{j-1}-2)|\nu|^{-1/j}
	\]
	holds in either of the following two cases:
	\begin{enumerate}[(i)]
		\item $j=1$ and $\Phi'(x)$ is monotone on $(a,b)$;
		\item $j\geq 2$.
	\end{enumerate}
\end{lem}

Let us turn to the proof of  Theorem~\ref{prop-uniform-main-esti-1}.

\begin{proof}[Proof of Theorem~\ref{prop-uniform-main-esti-1}]
	We first prove the uniform upper estimate.  Let $d,\eta$ and $r$ be fixed as in the statement and let us fix $B\in \mathcal{B}_{d,\eta,r}$. As in the proof of
	Theorem~\ref{Mainesti}~\eqref{Main-esti-1}, write
	\[
	B(e^{i\theta})=e^{i\psi_B(\theta)},\qquad \theta\in[0,2\pi].
	\]
	By continuity of  the map
	$(\Lambda,\theta)\mapsto \psi_{B_\Lambda}^{(q)}(\theta)$ on $K_{d,\eta,r}$,
	\[
	m_1:=\inf_{B\in\mathcal B_{d,\eta,r}}\inf_{\theta\in[0,2\pi]}\psi'_B(\theta)>0,\,\quad
	M_1:=\sup_{B\in\mathcal B_{d,\eta,r}}\sup_{\theta\in[0,2\pi]}\psi'_B(\theta)<\infty
	\]
	and
	\[
	M_2:=\sup_{B\in\mathcal B_{d,\eta,r}}\sup_{\theta\in[0,2\pi]}|\psi_B''(\theta)|<\infty.
	\]
	For \(n\ge1\) and \(k\ge0\), write
	\[
	\widehat{B^n}(k)
	=\frac1{2\pi}\int_0^{2\pi}e^{in\Phi_{B,a}(\theta)}\,d\theta,
	\quad
	\Phi_{B,a}(\theta):=\psi_B(\theta)-a\theta,
	\quad a:=\frac{k}{n}.
	\]
		Choose \(0<\delta<m_1/2\) and set
		\[
		K_\delta:=[m_1 -\delta,M_1+\delta].
		\]
		If $a<m_1-\delta$, then for every \(B\in\mathcal B_{d,\eta,r}\) and every \(\theta\in[0,2\pi]\),
		\[
		\Phi'_{B,a}(\theta)=\psi'_B(\theta)-a\ge m_1-a>\delta.
		\]
		If $a>M_1+\delta$, then for every \(B\in\mathcal B_{d,\eta,r}\) and every \(\theta\in[0,2\pi]\),
		\[
		\Phi'_{B,a}(\theta)=\psi'_B(\theta)-a\le M_1-a<-\delta.
		\]
		Therefore, if $a\notin K_\delta$, then 
		\[
		|\Phi'_{B,a}(\theta)|\ge\delta,
		\quad B\in\mathcal B_{d,\eta,r},\quad \theta\in[0,2\pi].
		\]
	
	An integration by parts then gives
	\[
	\int_0^{2\pi}e^{in\Phi_{B,a}(\theta)}\,d\theta
	=\Big[\frac{e^{in\Phi_{B,a}(\theta)}}{in\Phi'_{B,a}(\theta)}\Big]_{0}^{2\pi}
	-\frac1{in}\int_0^{2\pi}e^{in\Phi_{B,a}(\theta)}
	\Big(\frac1{\Phi'_{B,a}(\theta)}\Big)'\,d\theta.
	\]
	Since $\psi_B(2\pi)-\psi_B(0)=2\pi d$, $a=k/n$ and $\Phi_{B,a}'(2\pi)=\Phi_{B,a}'(0)$, the boundary term vanishes. Hence
	\[
	\left|\int_0^{2\pi}e^{in\Phi_{B,a}(\theta)}\,d\theta\right|
	\le \frac1n\int_0^{2\pi}
	\frac{|\Phi''_{B,a}(\theta)|}{|\Phi'_{B,a}(\theta)|^2}\,d\theta
	\le \frac{2\pi M_2}{n\delta^2}.
	\]
	Thus, for $a\notin K_\delta$, we obtain the stronger estimate
	\begin{equation}\label{ineq-outside-Kdelta}
		\left|\int_0^{2\pi}e^{in\Phi_{B,a}(\theta)}\,d\theta\right|\lesssim n^{-1},
	\end{equation}
	where the involved constants depend only on $d,\eta,r$, and not on $B$ or $a$.
	
	\medskip
	
	It remains to consider $a\in K_\delta$. By Lemma~\ref{lem:uniform-finite-type}, there exist
	$N_{d,\eta,r}\ge3$ and $c_{d,\eta,r}>0$, depending only on $d,\eta,r$, such that
	\begin{equation}\label{cetoile}
		\max_{2\le j\le N_{d,\eta,r}}|\psi_B^{(j)}(\theta)|\ge c_{d,\eta,r},
		\quad B\in\mathcal B_{d,\eta,r},\quad \theta\in[0,2\pi].
	\end{equation}
	We also have
	\[
	M:=\max_{2\le j\le N_{d,\eta,r}+1}
	\sup_{B\in\mathcal B_{d,\eta,r}}\sup_{\theta\in[0,2\pi]}|\psi_B^{(j)}(\theta)|<\infty.
	\]
	Choose \(\ell>0\) so small that \(M\ell\le c_{d,\eta,r}/2\), and partition \([0,2\pi]\) by
	finitely many intervals \(J_1,\ldots,J_Q\) of length at most \(\ell\).  Note that the number
	\(Q\) depends only on \(d,\eta,r\).  Fix one such interval \(J_q\) and a point
	\(\theta_q\in J_q\).  By \eqref{cetoile}, there is
	\(j_q\in\{2,\ldots,N_{d,\eta,r}\}\) such that
	\(|\psi_B^{(j_q)}(\theta_q)|\ge c_{d,\eta,r}\).  The mean-value theorem and the definition of $M$ give, for \(\theta\in J_q\),
	\[
	|\Phi_{B,a}^{(j_q)}(\theta)|
	=|\psi_B^{(j_q)}(\theta)| \geq |\psi_B^{(j_q)}(\theta_q)| - M\ell
	\ge \frac{c_{d,\eta,r}}{2}.
	\]
	Lemma \ref{lem:vdc}~(ii) applied on each interval $J_q$ yields
	\[
	\left|\int_{J_q}e^{in\Phi_{B,a}(\theta)}\,d\theta\right|
	\le 5\cdot2^{N_{d,\eta,r}}(c_{d,\eta,r}n)^{-1/j_q}.
	\]
	After summing over the $Q$ intervals $J_q$, we obtain, for any $a=k/n\in K_{\delta}$,
	\[
	|\widehat{B^n}(k)|\le \frac{1}{2\pi}\sum_{q=1}^Q\left|\int_{J_q}e^{in\Phi_{B,a}(\theta)}\,d\theta\right| \leq C_{d,\eta,r}n^{-1/N_{d,\eta,r}},
	\]
	where $C_{d,\eta,r}$ depends only on $d,\eta$ and $r$. Together with \eqref{ineq-outside-Kdelta}, this gives the stated upper estimate.
	
	\medskip
	
	The proof of the tail estimates follows from Cauchy's inequalities. Let $0<\rho<1<R<1/r$ and set
	\[
	m_\rho:=\sup_{B\in\mathcal B_{d,\eta,r}}\max_{|z|=\rho}|B(z)|<1,\qquad
	M_R:=\sup_{B\in\mathcal B_{d,\eta,r}}\max_{|z|=R}|B(z)|<\infty.
	\]
		Choose $0<J<1<L$ so that $m_\rho<\rho^J$ and $M_R<R^L$. Cauchy's estimates give, uniformly in $B\in\mathcal B_{d,\eta,r}$,
		\[
		\sum_{k<Jn}|\widehat{B^n}(k)|^2\le
		\frac{1}{1-\rho^2}\left(\frac{m_\rho}{\rho^J}\right)^{2n},
		\qquad
		\sum_{k>Ln}|\widehat{B^n}(k)|^2\le
		\frac{R^2}{R^2-1}\left(\frac{M_R}{R^L}\right)^{2n}.
		\]
		This ensures the existence of the desired $q\in (0,1)$ and $n_0\in \N$.
\end{proof}

\section{A uniform Hal\'asz type result for finite Blaschke products}\label{sec-Halasz}

We start with a general easy result. We recall that a subset of a topological space $X$ is a dense $G_\delta$ set if it is a countable intersection of dense open subsets of $X$.

\begin{prop}\label{prop-general-Cphi-Halasz}Let $p,q\in [1,\infty]$, let $w$ be a weight and let $\varphi:\mathbb{D}\to \mathbb{D}$ be holomorphic. If $C_{\varphi}$ does not map $\ell_a^p(w)$ into $\ell_a^q(w)$, then the set
	\[
	\UU:=\{f\in \ell_a^p(w):\,f\circ \varphi \notin \ell_a^q(w)\}
	\]
	is a dense $G_{\delta}$ subset of $\ell_a^p(w)$.
\end{prop}

\begin{proof}We assume that $p,q\in[1,\infty)$, the proof being similar for $p=\infty$ or $q=\infty$. For $N\geq 0$, define the linear operator $T_N:\ell_a^p(w)\to\ell^q$ by
\[
T_N(f):=\big(w_k\widehat{f\circ\varphi}(k)\mathbf1_{\{k\le N\}}\big)_{k\ge0}.
\]
Each $T_N$ is continuous. Indeed, for $0<r<1$, put $M_r:=\max_{|z|=r}|\varphi(z)|<1$. Cauchy's inequalities give
\[
|\widehat{\varphi^n}(k)|\le M_r^n r^{-k},\quad n,k\ge0.
\]
Since $\liminf_n w_n^{1/n}\ge1$, the sequence $(\widehat{\varphi^n}(k)/w_n)_n$ belongs to $\ell^{p'}$ for each fixed $k$, hence the linear map
\[
f\mapsto\widehat{f\circ\varphi}(k)=\sum_{n\geq 0}\widehat{f}(n)\widehat{\varphi^n}(k)
\]
is continuous on $\ell_a^p(w)$, by Hölder's inequality.

By assumption, there exists $f_0\in\ell_a^p(w)$ such that $f_0\circ\varphi\notin\ell_a^q(w)$. Consequently,
\[
\sup_{N\ge0}\|T_N(f_0)\|_{\ell^q}=\infty
\]
and thus $\sup_{N\ge0}\|T_N\|=\infty$. Now, for $m\geq 1$, set
\[
E_m:=\{f\in\ell_a^p(w):\ \sup_{N\ge0}\|T_N(f)\|_{\ell^q}\le m\}.
\]
 Note that
 \[
 \UU=\bigcap_{m\geq 1}E_m^c
 \]
 hence, by the Baire Category Theorem, it suffices to prove that each $E_m$ is closed and has empty interior. The fact that $E_m$ is closed follows from the continuity of $T_N$, $N \geq 0$. If $E_m$ contained a ball of radius $\varepsilon>0$, linearity would imply $\sup_N\|T_N\|\le2m/\varepsilon$, a contradiction.
\end{proof}

We recall that a weight $w$ is called slowly oscillating if for every $0<J<1<L$ there exists constants $c,C>0$ such that, for every $n$ large enough and every $Jn\leq k\leq Ln$, one has
\[
cw_n \leq w_k \leq Cw_n.
\]
It is not difficult to check that slowly oscillating weights have polynomial bound, from above and from below. More precisely,
\begin{prop}[Proposition~2.5 in \cite{LLQR-l2weight-2025}]\label{polyboundweight}Let $w$ be a slowly oscillating weight. Then there exist $\alpha ,\beta >0$ such that, for every $n$ large enough,
\[	
n^{-\alpha} \leq w_n \leq n^{\beta}.
\]
\end{prop}

Combining the previous with Baire's theorem and the fact that, given any non-monomial finite Blaschke product $B$, $C_B$ does not map $\ell_a^1(w)$ into $\ell_a^p(w)$, $p\in[1,2)$, nor $\ell_a^p(w)$ into $\ell_a^\infty(w)$ (as recalled in the introduction), we obtain:

\begin{cor}\label{cor-generic-unboundedness}Let $w$ be a slowly oscillating weight and let $B$ be a non-monomial finite Blaschke product.
	\begin{enumerate}[(a)]
		\item For every $p\geq 1$, there exists a dense $G_{\delta}$-subset $\mathcal{G}_B$ of $\ell_a^p(w)$ such that for every $f \in \mathcal{G}_B$, one has
		\[
		f\circ B \notin \bigcup_{q\in [1,2)}\ell_a^q(w).
		\]
		\item For every $p>2$, there exists a dense $G_{\delta}$-subset $\mathcal{H}_B$ of $\ell_a^p(w)$ such that for every $f \in \mathcal{H}_B$, one has
		\[
		f\circ B \notin \bigcup_{q\geq 1}\ell_a^q(w).
		\]
	\end{enumerate}	
	
\end{cor}

The previous application of Baire's theorem naturally raises the question of whether the sets $\mathcal{G}_B$ or $\mathcal{H}_B$  may be chosen independently of the finite Blaschke product $B$. Theorem~A, restated below, gives a positive answer for $\GG_B$. The corresponding question for $\HH_B$ remains open. Besides extending Hal\'asz's theorem (Theorem~\ref{Halaszthm}) to arbitrary non-monomial finite Blaschke products, it also provides a quantitative strengthening of Corollary~\ref{cor-generic-unboundedness}.

\begin{thm}\label{thm-Halasz-generic}
	Let $p\geq 1$, let $w$ be a slowly oscillating weight, and let $(v_k)_{k\geq 0}$ be a sequence of positive real numbers such that $v_k\to \infty$ as $k\to\infty$. Then, there exists a dense $G_{\delta}$-subset $\mathcal{G}$ of $\ell_a^p(w)$ such that, for every $f\in \mathcal{G}$ and every non-monomial finite Blaschke product $B$, one has
	\begin{equation}\label{halaszeq}
	\sum_{k\geq 1}|\widehat{f\circ B}(k)w_k|^{2-\frac{v_k}{\log (k+1)}}=\infty.
	\end{equation}
	In particular, there exists a dense $G_\delta$-subset $\tilde{\GG}$ of $\ell_a^p(w)$ such that for every $f\in \tilde{\GG}$ and every non-monomial finite Blaschke product $B$,
	\[
	f\circ B \notin \bigcup_{q\in [1,2)} \ell_a^q(w).
	\]
\end{thm}

It is not difficult to see that, if $(v_k)_{k\geq 0}$ is bounded, then the result fails in general (see \cite{Halasz1967}).

\begin{proof}[Proof of Theorem~\ref{thm-Halasz-generic}]
Fix $p\ge1$ and a slowly oscillating weight $w$. By the Baire Category Theorem, it is enough to establish the result uniformly for every Blaschke product in each compact family $K:=\mathcal B_{d,\eta,r}$, $d\geq 1$ and $0<\eta<r<1$ rational. We write
\[
\lambda_k:=2-\frac{v_k}{\log(k+1)},\quad k\geq 1.
\]
For $m\geq 1$, set
\[
E_m:=\left\{f\in\ell_a^p(w):\,\exists B\in K,\,\sum_{k\geq 1}|w_k \widehat{f\circ B}(k)|^{\lambda_k}\leq m\right\}.
\]
Observe that, if we define
\[
\GG=\ell_a^p(w) \setminus \bigcup _{m\geq 1}E_m,
\]
then every $f \in \GG$ satisfies \eqref{halaszeq}. Thus, by using Baire's theorem again, it suffices to prove that $E_m$ is closed and has empty interior.

Proceeding as in the proof of Proposition~\ref{prop-general-Cphi-Halasz}, it is easily proved that, for each fixed $k$, the map $(B,f)\mapsto \widehat{f\circ B}(k)$ is continuous on $K\times\ell_a^p(w)$. Moreover, with the above conventions, the functions $t\mapsto t^{\lambda_k}$ are continuous from $[0,\infty)$ to $[0,\infty]$. Let $f_j\in E_m$ converge to $f$ in $\ell_a^p(w)$, and let $B_j\in K$ be such that
\[
\sum_{k\geq 1}|w_k \widehat{f_j\circ B_j}(k)|^{\lambda_k}\leq m.
\]
By compactness and upon extracting, we may and shall assume that $B_j\to B\in K$. Thus
\[
\sum_{k\ge1}|w_k \widehat{f\circ B}(k)|^{\lambda_k} =  \sum_{k\ge1}\liminf_{j\to\infty}|w_k \widehat{f_j\circ B_j}(k)|^{\lambda_k}\leq \liminf_{j\to\infty}\sum_{k\ge1}|w_k \widehat{f_j\circ B_j}(k)|^{\lambda_k}\leq m,
\]
which proves that $f\in E_m$ and that $E_m$ is closed.

We now prove that each $E_m$ has empty interior. First, observe that this follows easily if $\lambda_k\leq0$ for infinitely many $k$. Indeed, for every polynomial $P$ and every $B\in K$, by Lemma~\ref{lem-easy-expo-decay}, $|\widehat{P\circ B}(k)|\leq CR^{-k}$ for some $R>1$ and $C>0$ uniform in $B\in K$. Since $w$ has at most polynomial growth by Proposition~\ref{polyboundweight}, we thus have
\[
|w_k\widehat{P\circ B}(k)|^{\lambda_k} \geq 1,
\]
for all sufficiently large $k$ such that $\lambda_k\leq0$. Hence $P\notin E_m$ for every $m\geq1$ and, since polynomials are dense in $\ell_a^p(w)$, the closed set $E_m$ has empty interior.

In the rest of the proof, we may assume that $\lambda_k>0$ for all sufficiently large $k$. By Theorem~\ref{prop-uniform-main-esti-1} and Parseval's identity, there exist $0<J<1<L$, $N\ge3$ and $C>0$ such that, for every $n$ sufficiently large and every $B\in K$,
\begin{equation}\label{eqthm330}
\sup_{k\geq 0}|\widehat{B^n}(k)|\leq Cn^{-1/N} \quad \text{and} \quad \sum_{Jn\leq k\leq Ln}|\widehat{B^n}(k)|^2\geq\frac12.
\end{equation}
By the first inequality in \eqref{eqthm330}, and since $Cn^{-1/N}\leq n^{-1/(2N)}$ and $\log(Ln+1)\leq2\log n$ for $n$ large enough, we get
\begin{equation}\label{eqthm33II}
	|\widehat{B^n}(k)|^{\lambda_k} \geq |\widehat{B^n}(k)|^2 \left(\sup_{Jn\leq j\leq Ln}|\widehat{B^n}(j)|\right)^{-v_k/\log(k+1)} \geq |\widehat{B^n}(k)|^2\exp\left(\frac{1}{4N}\min_{Jn\leq j\leq Ln}v_j\right),
\end{equation}
whenever $Jn\leq k\leq Ln$ and $n$ is large enough. Since the length of $[Jn,Ln]$ is $O(n)$, the second inequality in \eqref{eqthm330} implies, uniformly in $B\in K$ and for $n$ large enough,
\begin{equation}\label{eqthm33}
\sum_{\substack{Jn\leq k\leq Ln\\ |\widehat{B^n}(k)|\geq n^{-1}}} |\widehat{B^n}(k)|^2\geq\frac14.
\end{equation}
Thus, from \eqref{eqthm33II} and \eqref{eqthm33}, we deduce
\begin{equation}\label{eqthm33III}
\sum_{\substack{Jn\leq k\leq Ln\\ |\widehat{B^n}(k)|\geq n^{-1}}} |\widehat{B^n}(k)|^{\lambda_k} \geq\frac14\exp\left(\frac{1}{4N}\min_{Jn\leq k\leq Ln}v_k\right).
\end{equation}

Suppose now by contradiction that $E_m$ has nonempty interior. Since polynomials are dense in $\ell_a^p(w)$, there exist a polynomial $P$ and $0<\varepsilon< 1$ such that, for every $n\geq 0$,
\[
g_n(z):=P(z)+\frac{\varepsilon}{2w_n}z^n\in E_m.
\]
Since there exists $R>1$ such that every $B\in K$ is holomorphic on $|z|\le R$, by Lemma~\ref{lem-easy-expo-decay} and Proposition~\ref{polyboundweight}, we get, for every $n$ large enough (independent of $B\in K$),
\[
|\widehat{P\circ B}(k)| \leq \frac{\varepsilon}{4nw_n} \leq \frac{\varepsilon|\widehat{B^n}(k)|}{4w_n},
\] 
 for every $Jn\leq k\leq Ln$ such that $|\widehat{B^n}(k)|\geq n^{-1}$. Therefore, for all sufficiently large $n$, uniformly in $B\in K$ and for indices $Jn\leq k\leq Ln$ satisfying $|\widehat{B^n}(k)|\geq n^{-1}$, we get
\[
w_k|\widehat{g_n\circ B}(k)| \geq\frac{\varepsilon}{4}\frac{w_k}{w_n}|\widehat{B^n}(k)|  \gtrsim |\widehat{B^n}(k)|.
\]
Since $0<\lambda_k<2$ for $k$ large, the inequalities \eqref{eqthm33} and \eqref{eqthm33III} yield, for $n$ sufficiently large uniformly in $B\in K$,
\[
\begin{aligned}
\sum_{k\geq1}|w_k\widehat{g_n\circ B}(k)|^{\lambda_k}
& \gtrsim \sum_{\substack{Jn\leq k\leq Ln\\ |\widehat{B^n}(k)|\geq n^{-1}}} |\widehat{B^n}(k)|^{\lambda_k}\\
& \gtrsim \exp\!\left(\frac{1}{4N}\min_{Jn\leq k\leq Ln}v_k\right).
\end{aligned}
\]
The right-hand side goes to $+\infty$, since $v_k\to+\infty$. Thus, $g_n\notin E_m$ for $n$ large, a contradiction. This concludes the proof.
\end{proof}

\begin{rem}{\rm We point out that the divergence \eqref{halaszeq} in Theorem~\ref{thm-Halasz-generic} is actually uniform on compact families of Blaschke products. More precisely, for every $f\in\mathcal G$, $d\geq1$ and $0<\eta<r<1$, one has 
	\[
	\lim_{N\to\infty} 	\inf_{B\in\mathcal B_{d,\eta,r}} \sum_{k=1}^N |w_k\widehat{f\circ B}(k)|^{2-\frac{v_k}{\log(k+1)}}=+\infty.
	\]
	Indeed, for every $M>0$, the sets
	\[
	\left\{B\in\mathcal B_{d,\eta,r}:\sum_{k=1}^N|w_k\widehat{f\circ B}(k)|^{2-\frac{v_k}{\log(k+1)}}\leq M\right\}
	\]
	are closed and decreasing as $N$ increases, and their intersection is empty by Theorem~\ref{thm-Halasz-generic}. By compactness, one of these sets must be empty.
}
\end{rem}

\section{A $\ell_a^p(w)\to \ell_a^{\infty}(w)$ Hal\'asz type result for a fixed finite Blaschke product}\label{sec-dual-fixed-B}

This section is devoted to the proof of Theorem~B. We shall use the following weighted estimate on the columns of the matrix of a composition operator induced by a finite Blaschke product.  It follows from the column estimates established in \cite[Theorem~3.4]{CharpentierZarouf} and Cauchy's inequalities. For $p\in [1,\infty)$, we shall write $p'$ the conjugate exponent of $p$ (i.e. $1/p+1/p'=1$).

\begin{lem}\label{lem-column-growth-fixed-B}Let $w$ be a slowly oscillating weight, let $p>2$, and let $B$ be a non-monomial finite Blaschke product. Then there exists a constant $c>0$ and an increasing sequence $(k_j)_j$ tending to $\infty$ such that
\begin{equation}\label{eq-column-lower-fixed-B}
\left(\sum_{n\geq 0}|\widehat{B^n}(k_j)|^{p'}\frac{w_{k_j}^{p'}}{w_n^{p'}}\right)^{1/p'}\geq c k_j^{\frac12 -\frac1p},\quad j\geq 1.
\end{equation}
\end{lem}

\begin{proof}By \cite[Theorem~3.4(2)]{CharpentierZarouf}, there exist $0<J<L<\infty$ and $c>0$ such that, for every sufficiently large $N$,
\[
\frac1N\sum_{JN\leq k\leq LN}\sum_{n\geq 0}|\widehat{B^n}(k)|^{p'}\geq cN^{1-p'/2}.
\]
Hence, for every sufficiently large $N$, there exists $k_N\in[JN,LN]$ such that
\[
\left(\sum_{n\geq 0}|\widehat{B^n}(k_N)|^{p'}\right)^{1/p'}\gtrsim N^{\frac{1}{p'}-\frac12}\gtrsim k_N^{\frac12-\frac1p}.
\]
By Lemma~\ref{lem-elementary}, there exist $0<J<1<L$ such that, for every $\gamma >0$, there exists $0<s<1$ such that, for $k$ large enough,
\[
\sum_{n<Jk}\left|\widehat{B^{n}}(k)\right|^{p'} n^\gamma + \sum_{n>Lk}\left|\widehat{B^{n}}(k)\right|^{p'}n^\gamma \leq s^k.
\]
This, together with the fact that $w$ is slowly oscillating and thus admit polynomial bound below and above, we obtain, for every $N$ large enough,
\begin{equation*}
\sum_{n\geq 0}|\widehat{B^n}(k_N)|^{p'}\frac{w_{k_N}^{p'}}{w_n^{p'}} \geq \sum_{Jk_N\leq n\leq Lk_N}|\widehat{B^n}(k_N)|^{p'}\frac{w_{k_N}^{p'}}{w_n^{p'}} \gtrsim  \sum_{Jk_N\leq n\leq Lk_N}|\widehat{B^n}(k_N)|^{p'}\gtrsim \sum_{n\geq 0}|\widehat{B^n}(k_N)|^{p'}.
\end{equation*}

The result follows upon taking an increasing subsequence of $(k_N)_{N\geq 0}$.
\end{proof}

We are now ready to prove Theorem~B, that we restate as follows.

\begin{thm}\label{thm-dual-fixed-B}
Let $B$ be a non-monomial finite Blaschke product, let $p>2$, let $w$ be a slowly oscillating weight and let $(v_k)_{k\geq 1}$ be a sequence of positive real numbers such that
\begin{equation}\label{eq-subpolynomial-v}
v_k\longrightarrow\infty
\quad\text{and}\quad
v_k=o(k^{\frac{1}{2}-\frac{1}{p}}),\quad \text{as }k\to \infty.
\end{equation}
Then there exists a dense $G_\delta$-subset $\HH_{B,p,w}$ of $\ell_a^p(w)$ such that, for every $f\in\HH_{B,p,w}$,
\begin{equation}\label{eq-dual-fixed-B-main}
\limsup_{k\to\infty}\frac{|\widehat{f\circ B}(k)|w_k}{v_k}=\infty.
\end{equation}
In particular, for every $f\in\HH_{B,p,w}$, one has $f\circ B\notin\ell_a^\infty(w)$.
\end{thm}

\begin{proof}
For $k\ge1$, consider the linear functional $T_k:\ell_a^p(w) \to \C$ defined by
\[
T_k(f):=\frac{w_k}{v_k}\widehat{f\circ B}(k) = \frac{w_k}{v_k}\sum_{n\geq 0}\widehat{f}(n)\widehat{B^n}(k).
\]
Note that each $T_k$ is bounded on $\ell_a^p(w)$. Indeed, for $0<\rho<1$, put $M_\rho=\max_{|z|=\rho}|B(z)|<1$. Cauchy's inequality gives
\[
|\widehat{B^n}(k)|\le M_\rho^n\rho^{-k},
\]
hence, using that $w$ is slowly oscillating and Proposition~\ref{polyboundweight}, $(\widehat{B^n}(k)/w_n)_{n\geq 0}\in\ell^{p'}$ for every fixed $k$. The boundedness of $T_k$ then follows from Hölder's inequality.  Since the dual space of $\ell_a^p(w)$ identifies with $\ell_a^{p'}(w^{-1})$ for the Cauchy pairing that comes from the usual inner product of $\ell_a^2$, we have
\[
\Vert T_k \Vert = \frac{w_k}{v_k}\left(\sum_{n\geq 0}|\widehat{B^n}(k)|^{p'}w_n^{-p'}\right)^{1/p'}.
\]
Now, by Lemma~\ref{lem-column-growth-fixed-B}, there is a sequence $k_j\to\infty$ such that
\[
\Vert T_{k_j}\Vert=\frac1{v_{k_j}}\left(\sum_{n\ge0}|\widehat{B^n}(k_j)|^{p'} \frac{w_{k_j}^{p'}}{w_n^{p'}} \right)^{1/p'} \geq c\,\frac{k_j^{\frac12-\frac1p}}{v_{k_j}}.
\]
By the assumption \eqref{eq-subpolynomial-v} on $(v_k)_{k\geq 0}$, the right-hand side goes to $+\infty$ as $j\to \infty$.

Now, for $m\ge1$ let
\begin{equation}\label{eq-def-Em-dual}
E_m:=\{f\in\ell_a^p(w):\ |T_k(f)|\le m\text{ for every }k\geq 1\}.
\end{equation}
Each $E_m$ is closed by continuity of $T_k$. If $E_m$ contained a nonempty ball of radius $\varepsilon$, then, by linearity, it would imply $\sup_k\|T_k\|\le2m/\varepsilon$, contradicting the previous estimate. Thus, the set
\[
\HH_{B,p,w}:=\ell_a^p(w)\setminus\bigcup_{m\ge1}E_m
\]
is a dense $G_\delta$ subset of $\ell_a^p(w)$. To conclude, observe that every $f\in\HH_{B,p,w}$ satisfies \eqref{eq-dual-fixed-B-main}.
\end{proof}

We shall now derive a similar result, replacing $\ell_a^p(w)$, $p>2$, with the Fr\'echet space
\[
X(w):=\bigcap_{p>2}\ell_a^p(w),
\]
endowed with the topology induced by the increasing family of norms $\Vert \cdot \Vert_{\ell_a^{2+1/j}(w)}$. Recall that the polynomials are dense in $X(w)$.

\begin{cor}\label{cor-dual-frechet-weighted}
Let $w$ be a slowly oscillating weight and let $(v_k)_{k\geq 0}$ be a sequence of positive numbers. Assume that
\begin{equation}\label{eq-subpolynomial-v-inter}
	v_k\longrightarrow\infty \quad\text{and}\quad\frac{\log v_k}{\log(k+1)}\longrightarrow 0.
\end{equation}
Then, for every non-monomial finite Blaschke product $B$, there exists a dense $G_\delta$-subset $\HH_{B,w}$ of $X(w)$ such that, for every $f\in\HH_{B,w}$,
\[
\limsup_{k\to\infty}\frac{|\widehat{f\circ B}(k)|w_k}{v_k}=\infty.
\]
\end{cor}

\begin{proof}Note that \eqref{eq-subpolynomial-v-inter} implies \eqref{eq-subpolynomial-v} in Theorem~\ref{thm-dual-fixed-B}. Write $X_j=\ell_a^{2+1/j}(w)$, so that $X(w)=\bigcap_{j\geq1}X_j$. For $m\geq1$, define
	\[
	E_m:=\{f\in X(w):\ |T_k(f)|\leq m\text{ for every }k\geq1\},
	\]
	where $T_k$ is as in the proof of Theorem~\ref{thm-dual-fixed-B}. The sets $E_m$ are closed in $X(w)$. Observe that $E_m=X(w)\cap E_m^{(j)}$ for every $j$, where
	\[
	E_m^{(j)}:=\{f\in X_j:\ |T_k(f)|\leq m\text{ for every }k\geq 1\}.
	\]
	By Theorem~\ref{thm-dual-fixed-B}, applied in each $X_j$ and since \eqref{eq-subpolynomial-v} holds, the set $E_m^{(j)}$ is closed with empty interior in $X_j$. Suppose that $E_m$ has nonempty interior in $X(w)$. Since the norms $\Vert \cdot \Vert_{X_j}$ are increasing, there exist $j\geq1$, $f_0\in X(w)$ and $\varepsilon>0$ such that
	\[
	\{f\in X(w):\ \|f-f_0\|_{X_j}<\varepsilon\}\subset E_m.
	\]
	But $X_j\setminus E_m^{(j)}$ is open and dense in $X_j$, and $X(w)$ is dense in $X_j$ (it contains the polynomials). Therefore, there exists $f\in X(w)\setminus E_m^{(j)}$ with $\|f-f_0\|_{X_j}<\varepsilon$, a contradiction. Thus every $E_m$ is nowhere dense in $X(w)$, and the Baire Category Theorem completes the proof.
\end{proof}

\begin{rem}\label{rem-optimal-subpolynomial}{\rm The condition $\log(v_k)/\log(k+1)\to 0$ as $k\to \infty$ in Corollary~\ref{cor-dual-frechet-weighted} is somehow optimal. Indeed, let $f\in X(w)$ and let $\delta>0$. Choose $p>2$ so close to $2$ that
\[
\frac12 -\frac1p<\delta.
\]
Proceeding as in the proof of  Lemma~\ref{lem-column-growth-fixed-B}, one can easily derive from \cite[Theorem~3.4(1)]{CharpentierZarouf} the estimate
\begin{equation}\label{eq-column-upper-fixed-B}
	\left(\sum_{n\geq 0}|\widehat{B^n}(k)|^{p'}\frac{w_k^{p'}}{w_n^{p'}}\right)^{1/p'}\leq C(k+1)^{\frac12 -\frac1p}.
\end{equation}
Then, together with H\"older's inequality, \eqref{eq-column-upper-fixed-B} yields
\[
w_k|\widehat{f\circ B}(k)|\leq \|f\|_{\ell_a^p(w)}\left(\sum_{n\geq 0}|\widehat{B^n}(k)|^{p'}\frac{w_k^{p'}}{w_n^{p'}}\right)^{1/p'}\leq C(k+1)^{\frac12 -\frac1p}=o(k^\delta).
\]
Thus every $f\in X(w)$ satisfies $w_k|\widehat{f\circ B}(k)|=O(k^\delta)$ for every $\delta>0$. In particular, the conclusion of Corollary~\ref{cor-dual-frechet-weighted} does not hold in general, if one replaces $\log(v_k)/\log(k+1)\to 0$ by $v_k/k^a\to 0$ for any $a>0$.
}
\end{rem}

\medskip

Let us state the remaining open question.

\begin{question}\label{question-common-single-factors}
Let $(v_k)_{k\geq 0}$ satisfy \eqref{eq-subpolynomial-v-inter}. Does there exist a dense $G_\delta$ subset $\HH$ of $\cap_{p>2}\ell_a^p$ such that, for every $f\in\HH$ and every $\lambda\in\D\setminus\{0\}$,
\[
\limsup_{k\to\infty}\frac{|\widehat{f\circ b_\lambda}(k)|}{v_k}=\infty?
\]
The question remains open even for the weaker requirement corresponding to $v_k=1$, $k\geq 0$.
\end{question}

\bibliographystyle{amsplain}
\bibliography{refs-new}

\end{document}